\documentclass[12pt]{article}

\usepackage[english]{babel}
\usepackage{amssymb}
\usepackage{amsfonts}
\usepackage{amsmath}
\usepackage{bbold}

\usepackage{cite}

\renewcommand{\baselinestretch}{1.3}
\makeatletter
\def\@begintheorem#1#2{\trivlist%
 \item[\hskip \labelsep{\sffamily\bfseries #2\ #1}]\itshape}

\newtheorem{teo}{Theorem}[section]
\newtheorem{defi}[teo]{Definition}

\newtheorem{lem}[teo]{Lemma}
\newtheorem{pro}[teo]{Proposition}
\newtheorem{_rem}[teo]{Remark}
\newtheorem{_rems}[teo]{Remarks}
\newtheorem{_eje}[teo]{Example}

\newenvironment{rem}{\def\@begintheorem##1##2{\trivlist%
 \item[\hskip\labelsep{\sffamily\bfseries ##2\ ##1}]}\begin{_rem}}{\end{_rem}}
\newenvironment{rems}{\def\@begintheorem##1##2{\trivlist%
 \item[\hskip\labelsep{\sffamily\bfseries ##2\ ##1}]}
\begin{_rems}}{\end{_rems}}

\makeatother

\newenvironment{beweis}{{\noindent \em Proof:}}{\hfill $\rule{2mm}{3mm}$
\vspace{3mm}}

\DeclareMathAlphabet{\Ma}{U}{msa}{m}{n}
\DeclareMathAlphabet{\Mb}{U}{msb}{m}{n}
\DeclareMathAlphabet{\Meuf}{U}{euf}{m}{n}

\def\got#1{\Meuf{#1}}

\DeclareSymbolFont{ASMa}{U}{msa}{m}{n}
\DeclareSymbolFont{ASMb}{U}{msb}{m}{n}
\DeclareMathSymbol{\hrist}{\mathord}{ASMa}{"16}
\DeclareMathSymbol{\varkappa}{\mathalpha}{ASMb}{"7B}
\DeclareMathSymbol{\CrPr}{\mathord}{ASMb}{"6F}

\def\rest{\upharpoonright}

  \def\al #1.{{\mathcal{#1}}}
  \def\ot #1.{{\got{#1}}}
\def\CCRX{\al W.(X,\sigma)}
  \def\ccr #1,#2.{\overline{\Delta(#1,\,#2)}}

  \def\b #1.{{\bf #1}}
  \def\cross#1.{\mathrel{\mathop{\times}\limits_{#1}}}

  \def\C{\Mb{C}}
  \def\N{\Mb{N}}
  
  \def\R{\Mb{R}}
  \def\Ro{\Mb{R} \backslash \{0\}}

 \def\un{{\mbox{\boldmath $1$}}}

\def\f #1,#2.{\mathsurround=0pt \hbox{${#1\over #2}$}\mathsurround=5pt}

  \def\cross #1.{\mathrel{\raise 3pt\hbox{$\mathop\times\limits_{#1}$}}}
  \def\ol #1.{\overline{#1}}
\def\b #1.{{\bf #1}}

\def\ker{{\rm Ker}\,}
\def\Her #1.{{\rm Her}(#1)}

\def\dom{{\rm Dom}\,}

\def\rlf{{R(\lambda,f)}}
\def\rsl{\mathord{\al R.(X,\sigma)}}
\def\s #1.{_{\smash{\lower2pt\hbox{\mathsurround=0pt 
$\scriptstyle #1$}}\mathsurround=5pt}}
\def\set #1,#2.{\left\{\,#1\;\bigm|\;#2\,\right\}}
\def\maprightu #1;{\smash{\mathop{\longrightarrow}\limits^{#1}}}
\def\maprightd #1;{\smash{\mathop{\longrightarrow}\limits_{#1}}}
\def\maprightt #1,#2.{\mathrel{\smash{\mathop{\longrightarrow}
\limits_{#1}^{#2}}}}

\def\j{\phi}

\def\XP#1!{\renewcommand{\baselinestretch}{.7}\marginpar
{$\leftarrow${\footnotesize #1}\hfil} \renewcommand{\baselinestretch}{1}}
\def\XB{\marginpar{
{\footnotesize Change~starts-}\lower 11pt\hbox{\mathsurround=0pt$
\!\!\displaystyle{
\Bigg\downarrow}$\mathsurround=3pt}}}
\def\XE{\marginpar{{\footnotesize Change~ends-}\raise 
10pt\hbox{\mathsurround=0pt$\!\!\displaystyle{\Bigg\uparrow}$
\mathsurround=3pt}}}
\def\f #1,#2.{\mathsurround=0pt \hbox{${#1\over #2}$}\mathsurround=5pt}

\def\affp #1.{\mathrel{\Finv\s{#1}.}}

\def\ie{\textit{i.e.\ }}
\def\eg{\textit{e.g.\ }}

\def\viz{\textit{viz.\ }}
\def\margin #1.{\marginpar{#1}}

\def\spec{\mbox{\rm spec}}

\title{\bf The Resolvent Algebra: Ideals and Dimension}
\author{{}\\[-5mm] {\Large Detlev Buchholz}\\[10mm]
Institut f\"ur Theoretische Physik and Courant Centre \\
``Higher Order Structures in Mathematics'',
Universit\"at G\"ottingen, \\ 37077 G\"ottingen, Germany}

\date{}

\begin{document}
\maketitle

\begin{abstract}
\noindent
Let  $(X,\sigma)$ be a symplectic space admitting a 
complex structure and let $\al R.(X,\sigma)$
be the corresponding resolvent algebra, \ie the 
C*--algebra generated by the resolvents of selfadjoint operators
satisfying canonical commutation relations associated with $(X,\sigma)$.
In previous work this algebra was shown to provide a convenient
framework for the analysis of quantum systems. In the present
article its mathematical properties are elaborated with 
emphasis on its ideal structure.  
It is shown that  $\al R.(X,\sigma)$ is always nuclear and, 
if $X$ is finite dimensional, also of type I (postliminal).
In the latter case $\dim(X)$ labels the isomorphism classes 
of the corresponding resolvent algebras. 
For $X$ of arbitrary dimension, 
principal ideals are identified which are the 
building blocks for all other ideals. The maximal and minimal
ideals of the resolvent algebra are also determined. \\[3mm]
{\bf Keywords:} \ 
resolvent algebra, symplectic space, 
canonical commutation relations, \mbox{C*--algebra}, 
ideal, representation, isomorphism class  
\end{abstract}

\vfill\eject

\section{Introduction}
\setcounter{equation}{0}

In \cite{BuGr2} we have defined and analyzed the 
{resolvent algebra} of the canonical
commutation relations. Apart from the applications in that paper, 
this algebra has already demonstrated
its usefulness elsewhere. For example, on its basis 
one can model in a C*--context superderivations which occur
in supersymmetry,  cf.~\cite{BuGr}, as well as in 
BRST-constraint theory, cf.~\cite{Cost}. It occurs also
naturally in the representation theory of abelian Lie
algebras of derivations acting on a C*–algebra \cite{BuGr3}. 
Here we continue our analysis of the resolvent
algebra, with particular emphasis on its ideal structure.

We review the background which motivates the 
study of the resolvent algebra.
Canonical systems of operators have always been a central ingredient
in the modeling of quantum systems.
These systems of operators may all be presented in the following
general form: there is
a real linear map $\j$  from a given symplectic space
${(X,\sigma)}$ to a linear space of selfadjoint operators
on some common dense invariant core $\al D.$ in a Hilbert space $\al H.$, 
satisfying the relations
$$
\big[\j(f),\,\j(g)\big]=i\sigma(f,\, g) \, {\un}, \quad 
\j(f)^* =\j(f)\quad\hbox{on}\quad
\al D. \, .
$$
In the case that $X$ is finite dimensional, one can reinterpret this
relation in terms of the familiar quantum mechanical
position and momentum operators, and if
$X$ consists of Schwartz functions on some 
manifold one may consider $\j$ to be a bosonic quantum field.
The observables of the system are then constructed from the operators
$\{\j(f) : f\in X\} \,,$ usually as polynomial expressions.
Since one wants to study a variety of representations of such
systems,  it is convenient to cast the algebraic
information of the canonical systems into C*--algebras,
given the rich source of mathematical tools available there.

The obvious way to take this step
is to form suitable bounded functions of the generically 
unbounded fields $\j(f)$.
In the approach introduced by Weyl, this is done by
considering the C*--algebra generated by the set of unitaries
$$
\big\{ \exp\big(i\j(f)\big) : f\in X\big\}
\, .
$$
Regarded as an abstract 
algebra, it is the familiar Weyl algebra  \cite{Man}, 
denoted by $\CCRX$. The Weyl algebra suffers, however, from 
several well--known flaws with regard to physics. First and foremost,
it does not admit the definition of much 
interesting dynamics (one--parameter automorphism groups),
cf.~\cite{FaVB,BuGr2}.  Second, natural observables
such as bounded functions of the Hamiltonian
are not in $\CCRX$. Third, the Weyl algebra has a 
vast number of representations 
in which representers of its generators $\j(f)$ cannot 
be defined \cite{Ace, Gr3, GrHu3}. These nonregular 
representations describe situations where the 
field $\j$ has ``infinite strength''. Whilst this is sometimes 
useful for idealizations, cf.\ for example the discussion
of plane waves in \cite{Ace} or of quantum constraints 
in \cite{GrHu3}, the majority of nonregular 
representations of the Weyl algebra is of no interest.

This motivates the consideration of alternative
C*--algebraic versions of the canonical commutation relations.
Instead of taking the C*--algebra generated
by exponentials of the underlying generators,
as for the Weyl algebra, we propose to consider 
the C*--algebra generated by their resolvents \cite{BuGr2}. 
These are given on the underlying Hilbert space $\al H.$ by
$$ \{ R(\lambda, f) \doteq (i \lambda \un - \j(f))^{-1} : 
\lambda \in \R\backslash \{ 0 \} \, , f \in X \} \, . $$
All algebraic properties of the fields can be expressed in terms of
relations amongst these resolvents and this fact
allows one to define the unital C*-algebra
generated by the resolvents also in representation independent terms. 
The structure of the resulting resolvent algebra will be studied below.

In contrast to the Weyl algebra, which is simple, the resolvent
algebra has ideals. This feature agrees with 
the observation that a unital C*--algebra which admits 
the definition of a sufficiently 
diverse variety of dynamics cannot be 
simple \cite[p.~2767]{BuGr2}. In the present article we analyze 
the ideal structure of the resolvent algebra and show 
that it depends sensitively on the size of the underlying quantum 
system. By a study of its primitive 
ideals we find that the resolvent algebra is  
of type I (postliminal) if the dimension of $X$ 
is finite and that it is merely nuclear if $X$ is infinite 
dimensional. Moreover, the specific nesting of its 
primitive ideals encodes information about the dimension of 
the underlying space $X$. As a matter of fact, this dimension,
if it is finite, is an algebraic invariant which labels the 
isomorphism classes of the associated resolvent algebras. We 
also analyze a distinguished family of principal ideals, generated by 
resolvents, which are the building blocks of all other ideals. 
Based on these results we determine the maximal and minimal
ideals of the resolvent algebras, the latter being  
trivial if $X$ is infinite dimensional. So, in summary, each resolvent 
algebra comprises very specific information 
about the underlying quantum system.

The article is organized as follows. We recall in the subsequent 
section some definitions and facts concerning 
the resolvent algebra which were established in \cite{BuGr2};
readers familiar with these results may skip this section. 
In Sect.~\ref{PrimitiveIdeals} we analyze the 
structure of its irreducible representations and of its 
primitive ideals. Sect.~\ref{PrincipalIdeals} contains the analysis of its  
elemental principal ideals and further resultant information
about its ideal structure. A brief discussion concludes this   
article.
 
\section{Resolvent algebra -- definitions and facts}
\setcounter{equation}{0}
\label{ResBasics}

For the convenience of the reader, 
we compile in this section some definitions and  
facts which were established in \cite{BuGr2}, where proofs
of the claims and further details can be found.
Let $X$ be a real vector space and let $\sigma : X \times X
\rightarrow \R$ be a nondegenerate symplectic form; in order to 
avoid pathologies, we make the standing assumption that 
$(X,\sigma)$ admits some complex structure \cite{Ro}.
Just as the Weyl algebra can be abstractly defined by the Weyl relations,
the C*-algebra of resolvents is abstractly defined 
by its generators and relations.
\begin{defi}
\label{ResAlg}
Given a symplectic space $(X,\,\sigma)$,  
$\al R._0$ is the universal unital
*-algebra generated by
the set $\{ R(\lambda,f) : \lambda \in \Ro, \; f\in X \}$
and the relations
\begin{eqnarray}
\label{Resolv}
\rlf - R(\mu,f) \! \! &=& \! \! i(\mu-\lambda)\rlf R(\mu,f)  \\[1mm]
\label{Rinvol}
R(\lambda,f)^*  \! \! &=& \! \! R(-\lambda,f) \\[1mm]
\label{Rccr}
\big[\rlf,\,R(\mu,g)\big] \! \!  &=& \! \! 
i\sigma(f,g)\,\rlf\,R(\mu,g)^2\rlf \\[1mm]
\label{Rhomog}
\nu \,R(\nu \lambda,\, \nu f)   \! \!  &=&  \! \!  \rlf \\[1mm]
\label{Rsum}
\rlf R(\mu,g) \! \! &=&  \! \! R(\lambda+\mu,\,f+g)\big(\rlf+R(\mu,g)
+i\sigma(f,g)\rlf^2R(\mu,g)\big) \\[1mm] 
\label{Riden}
R(\lambda,0) \! \! &=& \! \! -{\textstyle {i \over \lambda}}\,\un 
\end{eqnarray}
where $\lambda,\, \mu, \, \nu \in \Ro$ and $f,\,g\in X$,
and for (\ref{Rsum}) we require $\lambda+\mu\not=0$.
That is, start with the free unital *-algebra generated by
 $\{ R(\lambda,f) : \lambda \in \Ro, \; f \in X \}$
 and factor out by the ideal generated by the relations
(\ref{Resolv}) to (\ref{Riden}) to obtain the *-algebra $\al R._0$.
\end{defi}
\begin{rems}
(a) Relations (\ref{Resolv}), (\ref{Rinvol}) encode the 
algebraic properties of the resolvent of some selfadjoint operator,
(\ref{Rccr}) encodes the canonical commutation relations and relations 
(\ref{Rhomog}) to (\ref{Riden}) express the linearity of the underlying
map $\phi$ on $X$. \\ 
(b) \ The *-algebra $\al R._0$ is nontrivial, because it has nontrivial
representations. {}For instance, in a Fock representation
$\pi$ of the canonical commutation relations over $(X,\,\sigma)$
one has selfadjoint operators $\phi_\pi (f)$, $f \in X$, 
satisfying the commutation relations on a sufficiently big domain 
so that one can define $\pi(\rlf) \doteq (i\lambda\un-\phi_\pi (f))^{-1}$
to obtain a representation of $\al R._0.$ 
\end{rems}
Let $\ot S.$ denote the set of positive, normalized functionals 
$\omega$ of $\al R._0$. By \cite[Prop.~3.3]{BuGr2}  
their GNS--representations $(\pi_\omega, \al H._\omega)$ are uniformly 
bounded w.r.t.~$\ot S.$, so one can define: 
\begin{defi} \label{Rnorm}
Let $(X,\,\sigma)$ be a symplectic space and let $\al R._0$ be the 
corresponding *--algebra. Then
\[ \| A \| \doteq \sup_{\omega\in\ot S.} \, \|\pi_\omega(A)\|_{\al
  H._\omega} \, , \quad  A \in \al R._0
\] 
defines a C*--seminorm on $\al R._0$. The resolvent algebra 
$\rsl$  is defined as the C*--completion of the quotient 
algebra $\al R._0 \, / \, \ker \|\cdot\|$, where here and in the 
following the symbol $\, \ker \!$ denotes the kernel of 
the respective map. 
\end{defi}

\begin{rem}
It follows from \cite[Thm.~3.6]{BuGr2} that the functions 
$\lambda \mapsto R(\lambda,f)$ on $\R \backslash \{0\}$
can be analytically continued within $\rsl$
to the domain $\C \backslash i \R$. 
Obvious extensions of relations (\ref{Resolv}) to 
(\ref{Riden}) hold also for these continuations. 
\end{rem}

As mentioned, the resolvent algebra is not simple, in contrast 
to the Weyl algebra. In order to explore its rich ideal structure 
we need to study its representations. In this context the notion of a 
regular representation is of particular interest. 

\begin{defi}
A  representation $(\pi, \al H. \big)$ 
of $\rsl$ is regular on a set $S\subseteq X$ if
\[
\ker\pi\big(R(\lambda,f)\big)=\{0\} \quad \mbox{for all} \ f \in
S \ \mbox{and some (hence all)} \ \lambda \in \Ro \, .
\]
A state $\omega$ of $\rsl$ is {regular on} a set 
$S\subseteq X$ if its GNS--representation
$\pi_\omega$ is regular on $S\subseteq X$.
A representation (resp.\ state) which is regular on all of $X$
is said to be regular.
\end{defi}

This definition draws upon the fact that, 
as a consequence of the resolvent equations
(\ref{Resolv}),~(\ref{Rinvol}) 
all operators $\pi\big(R(\lambda,f)\big)$, $\lambda \in \R \backslash \{ 0 \}$,
have a common range and a common null space for fixed $f \in X$.
Moreover, if $\ker \pi\big(R(\lambda,f)\big) = \{ 0 \}$ for some (hence for
all) $\lambda \in \R \backslash \{ 0 \}$, they are resolvents
of some generator $\j_\pi(f)$, cf.~\cite[Ch.~VIII.4, Thm.~1]{Yos}. 
It is given by 
\begin{equation} \label{generator}
\j_\pi(f) \doteq i \lambda \un-\pi\big(R(\lambda,f)\big)^{-1} 
\end{equation} 
on the dense domain $\dom \j_\pi(f) = \pi\big(R(\lambda,f)\big) \, \al
H.$, $\lambda \in \Ro$. Some relevant properties of these 
(selfadjoint) operators have been established in \cite[Thm.~4.2]{BuGr2};
in particular, they provide a one--to--one correspondence  
between the regular representations of $\rsl$ and of~$\CCRX$. 
Further properties of the resolvent algebra $\rsl$ which are 
used in the present analysis are, cf.~\cite[Thm.~4.9]{BuGr2}:

\begin{pro}
\label{UniqueR}
Let $(X,\sigma)$ be a symplectic space of arbitrary dimension
and let $S \subset X$ be any finite dimensional 
nondegenerate subspace (including the zero--dimensional 
space $\{0\}$). Then:
\begin{itemize}
\item[(i)]
The norms of $\rsl$ and of $\al R.(S,\,\sigma)$ coincide on
${\hbox{*--Alg} \, \{\rlf : f\in S,\,\lambda\in\R\backslash \{0\}\}}$.
Thus one has the containment 
${\al R.(S,\,\sigma)}\subset \rsl.$
\item[(ii)]$\rsl$ is the inductive limit of the net of all
${\al R.(S,\,\sigma)}$ where $S\subset X$ ranges over all finite dimensional
nondegenerate subspaces of $X$ (including $X$ if it is finite dimensional).
\item[(iii)] Every regular representation of $\rsl$ is faithful.
\end{itemize}
\end{pro}

This result also has structural consequences. Since the 
irreducible Fock
representation induces a regular representation of $\rsl$, point (iii)
implies that $\rsl$ has faithful irreducible representations (\ie it is 
a primitive algebra). Thus its center must be trivial and 
every nonzero closed two--sided ideal $\al J. \subset \rsl$ is essential 
(\ie the left, hence also the right annihilator of $\al J.$ is zero  
since there exist cyclic vectors for $\al J.$
in every faithful irreducible representation). Consequently   
the intersection $\al J._1 \cap \al J._2$ of any two nonzero 
closed two--sided ideals $\al J._1, \, \al J._2 \subset \rsl$ is 
also nonzero (\ie $\rsl$ is prime), cf. Lemma \ref{InductIdeal}(i) below.

\section{Primitive ideals and the dimension of $X$}
\setcounter{equation}{0}
\label{PrimitiveIdeals}

In this section we determine the primitive ideals of $\rsl$ 
when the dimension $\dim(X)$ of $X$ is finite. 
The results enable us to show that $\dim(X)$ distinguishes the 
isomorphism classes of $\rsl$. 
As a further consequence we find 
that $\rsl$ is of type I (postliminal).

By definition, the primitive ideals are the kernels of 
irreducible representations, including the trivial ideal $\{0\}$. 
To determine them we rely on the subsequent lemma, 
established in \cite[Prop.~4.7]{BuGr2}, where we use following notation: 
The symplectic complement of any 
subspace $S\subseteq X$ is denoted
by $S^\perp\doteq \{ f\in X : \sigma(f,S)=0 \}$. 
The expression $X=S_1\oplus S_2\oplus\cdots\oplus S_n$ says that 
all subspaces $S_i \subset X$ are nondegenerate, 
$S_i\subseteq S_j^\perp$ if $i\not=j$, and each $f\in X$ has a unique 
decomposition
$f=f_1+f_2+\cdots+f_n$ with $f_i\in S_i$, $i = 1, \dots , n$.
Note that the zero--dimensional subspace~$\{0\} \subset X$ is admitted here for
notational convenience. 

\begin{lem}
\label{Xdecomp} 
\ Let $X$ be of finite or infinite dimension and 
let $(\pi,\al H.)$ be a representation of~$\rsl$.
\begin{itemize}
\item[(i)]
The set $X_R \doteq \{ f\in X : \ker\pi\big(\rlf\big)=\{0\} \}$  
is a linear space which is independent of 
$\lambda \in \Ro$. So for its complement 
$X_S\doteq X\backslash X_R$ one has $X_S + X_R = X_S$. 
\item[(ii)]
The set $X_T\doteq \{ f\in X : \ker\pi\big(\rlf\big)=\{0\}\;
\hbox{and}\;\pi\big(\rlf\big)^{-1}\in {\al B.(\al H.)} \} \subseteq
X_R$ is a linear space which is independent of $\lambda \in \Ro$. Moreover, 
$\pi\big(R(\mu,g)\big)=0$ for all $\mu \in \Ro$ and 
$g\in X$ with $\sigma(g,X_T)\not=0$.
Hence $\sigma(X_R,X_T)=0$.
\item[(iii)]
If $\pi$ is factorial, then $\pi\big(\rlf\big)=0$ for all $f\in X_S,$ and
$\pi\big(\rlf\big)\in\C \, \un\backslash \{0\}$ for all $f\in X_T$
and any $\lambda \in \Ro$. Moreover $X_T=X_R\cap {X_R}^\perp$.
\item[(iv)] Let $X$ be finite dimensional.
If $\pi$ is factorial one can augment $X_T$ by a
complementary space $X_T^{\, \sim} \subset X$ of the same dimension
such that the space $Q = X_T + X_T^{\, \sim}$ is nondegenerate. 
Moreover, one has the decomposition
\begin{equation}
\label{Qdecomp}
X=Q\oplus(Q^\perp\cap X_R)\oplus(Q^\perp\cap {X_R}^\perp)
\end{equation}
into nondegenerate spaces and 
$Q^\perp\cap X_R \subseteq (X_R\backslash X_T) \cup \{0\}$, \
$Q^\perp\cap {X_R}^\perp \subseteq X_S \cup \{0\}$.
\end{itemize}
\end{lem}
Clearly $X_R$ is the part of $X$ on which $\pi$ is regular,
$X_T$ is the part on which it is ``trivially regular'',
$X_S$ is the part on which it is singular, and these have a 
particularly nice form when $\pi$ is factorial. 
Refining this, we can in fact fully characterize all
irreducible representations and their respective kernels.
However, we first need to define an object which will be a convenient index set
of the unitary equivalence classes of irreducibles.

\begin{defi}
Given a symplectic space $(X,\sigma)$, 
let \ $I (X,\sigma)$ denote the set of all pairs 
${(Y, \, \chi)}$,
where $Y \subseteq X$ is any subspace and $\chi$ any pure state 
(character) on the abelian \mbox{C*--algebra} 
${C^*\{\rlf : \lambda \in \Ro  \, , \,  f \in Y \cap Y^\perp \}}$ 
which does not vanish
on any of its generating resolvents. Note that, as a consequence
of the relations (\ref{Resolv}) to (\ref{Riden}), there
is a linear functional $\varphi : Y \cap Y^\perp \rightarrow \R$
such that $\chi(R(\lambda, f)) = (i \lambda - \varphi(f))^{-1}$
for $\lambda \in \Ro$, $f \in Y \cap Y^\perp$. 
\end{defi}

Let $(\pi,\al H.)$ be any irreducible representation
of $\rsl$ and let $X_R, X_T$ be the associated subspaces
of $X$, defined in the preceding lemma. 
According to part (iii) of this lemma \mbox{$X_T=X_R\cap {X_R}^\perp$}
and $\pi({\al R.}_T) \in \C \, \un$, where 
$\al R._T \doteq 
{C^*\{\rlf : \lambda \in \Ro \, , \,  f\in X_T \}}$. 
Taking into account that $\pi$ is a homomorphism and 
$\al R._T$ an abelian C*--algebra 
it follows that there is a unique pure state 
$\chi_T$ on $\al R._T$ such that 
$\pi \rest {\al R._T} = \chi_T(\, \cdot \,) \, \un_{\al H.}$.  
Moreover, as an immediate consequence of the 
definition of $X_T$ and the resolvent equation (\ref{Resolv}), $\chi_T$
does not vanish on any of the generating resolvents of $\al R._T$. 
Thus one can assign to each irreducible 
representation of $\rsl$
a space and state with properties described in the
preceding definition:  
$(\pi,\al H.) \longmapsto (X_R,\chi_T) \in I (X,\sigma)$.
Note that all representations in the unitary equivalence class 
of $(\pi,\al H.)$
give rise to the same element of $I (X,\sigma)$
since the associated spaces $X_R,X_T$ 
and functionals $\chi_T$ remain fixed under 
arbitrary unitary transformations of the representation.
Thus, denoting the set of unitary equivalence classes 
of irreducible representations (the spectrum) of $\rsl$ by
$\widehat{\rsl}$, the preceding 
assignment determines a map
\[
\iota \, : \, \widehat{\rsl} \rightarrow I (X,\sigma) \, .
\] 
Its properties are studied in the subsequent proposition for
the case of finite dimensional symplectic spaces.

\begin{pro}
\label{Irrep}
Let $(X,\,\sigma)$ be any given symplectic space of finite dimension.
Then the map $\iota : \, \widehat{\rsl} \rightarrow 
I (X,\sigma)$ defined above is a bijection.
\end{pro}
\begin{beweis}
In a first step we prove that $\iota$ is injective. Let 
$(\pi, \al H.)$  and $(\pi^\prime, \al H.^\prime)$ be two 
irreducible representations of $\rsl$ determining the same 
element $(X_R,\chi_T) \in I (X,\sigma)$. We need
to show that these representations are equivalent.
Since a representation is uniquely specified by
its values for $\{R(\lambda,f) : \lambda \in \R\backslash \{0\} \, ,
f\in X\}$, we focus on  this set.
According to part (iii) of Lemma \ref{Xdecomp} and the resolvent equation
(\ref{Resolv}) 
we have $\pi(R(\lambda,f_S)) = 0$ for $f_S \in X_S = X\backslash X_R
\, ;$ moreover, $\pi(R(\lambda,f_T)) = \chi_T(R(\lambda,f_T)) \, \un_{\al H.}
= (i\lambda - \varphi(f_T))^{-1} \, \un_{\al H.}$ 
for  $f_T \in X_T = X_R \cap {X_R}^\perp$, where 
$\varphi : X_T \rightarrow \R$ is a linear functional. By
part (iv) of Lemma~\ref{Xdecomp}, there is 
a nondegenerate subspace $X_N \doteq Q^\perp \cap X_R \subseteq X_R$ 
(depending on the chosen extension Q of $X_T$) such that each $f_R \in X_R$
can uniquely be decomposed into $f_R = f_T + f_N$ where
$f_T \in X_T$ and $f_N \in X_N$. Plugging this into  
relation (\ref{generator}) and bearing in mind  the linearity of 
the generators $\phi_\pi$ we obtain
$\phi_\pi(f_R) = \varphi(f_T) \, \un_{\al H.} + \phi_\pi(f_N)$
which yields 
\begin{equation} \label{RepEq}
 \pi(R(\lambda, f_R)) = 
\pi(R(\lambda + i \varphi(f_T), f_N)) \, .
\end{equation}
It is an immediate consequence of these observations that
$\pi(\rsl) = \pi({\al R.} (X_N, \sigma)) $, hence
$\pi \rest {{\al R.} (X_N, \sigma)}$ is still irreducible. The same 
arguments apply to the representation $(\pi^\prime, {\al H.}^\prime)$, 
so one can exchange everywhere in the preceding relations
$\pi$ by $\pi^\prime$ and $\un_{\al H.}$ by $\un_{{\al H.}^\prime}$. 
Now since $X_N$ is finite dimensional and non--degenerate, all
regular irreducible representations of ${\al R.}(X_N,\sigma)$
are unitarily equivalent by the von Neumann uniqueness theorem,
cf. \cite[Cor.~4.4]{BuGr2}. Thus there is an isometry
$V$ mapping $\al H.$ onto ${\al H.}^\prime$ which intertwines
the representations $\pi \rest {{\al R.} (X_N, \sigma)}$ and
$\pi^\prime \rest {{\al R.} (X_N, \sigma)}$. The relations
established above then imply that 
$V \, \pi(R(\lambda, f)) = \pi^\prime (R(\lambda, f)) \, V$   
for all $\lambda \in \R\backslash \{0\}$, $f\in X$.
Hence the representations $(\pi, {\al H.})$ and 
$(\pi^\prime, {\al H.}^\prime)$ of the given algebra $\rsl$ are  
equivalent, proving the injectivity of $\iota$. 

For the proof of surjectivity, 
let  ${(Y,\chi)} \in I (X,\sigma)$ be given 
and let $Z \doteq Y \cap Y^\perp$. We have to show that there exists
some irreducible representation $(\pi,{\al H.})$ of $\rsl$ such that 
the associated spaces coincide with the given ones, \ie 
$X_R = Y$, $X_T = Z$, $X_S = X \backslash Y$, and 
$\pi \rest {\al R.}_T = \chi(\, \cdot \,) \, \un_{\al H.}$. 
In a first step we establish a decomposition of $X$ similar to 
that given in relation (\ref{Qdecomp}): Making use of 
\cite[Lem.~A.1(iii)]{BuGr2}, we pick some  
subspace $Z^{\, \sim} \subset X$ which has the same dimension as $Z$
and which is complementary to $Z$ in the sense that the 
linear space $Q \doteq Z + Z^{\, \sim}$ is nondegenerate.
Then each $f \in Y$ can uniquely be decomposed into 
$f = f_Z + f_N$ where $f_Z \in Z \subset Y$ and $ f_N 
\in Z^{\, \sim \perp} \cap Y = Q^\perp \cap Y$.
Moreover, the space $Q^\perp \cap Y$ is also nondegenerate; for 
if $\sigma(f_N,g) = \sigma(f_N,g_N) = 0$ 
for all $g = g_Z + g_N \in Y$ it follows that 
$f_N \in Q^\perp \cap Y \cap Y^\perp = Q^\perp \cap Z = \{0\}$.  

Now let $(\pi_0, {\al H.})$ be any regular irreducible 
representation of ${\al R.}(Q^\perp \cap Y, \sigma)$, \eg
the Schr\"odin\-ger representation, and let 
$\varphi : Z \rightarrow \R$ be the linear functional
fixed as above by the pure state~$\chi$ on 
$\al Z. \doteq {C^*\{\rlf :  \lambda \in \Ro \, , \, f \in Z \}}$. With 
this input and the preceding information we can define 
a representation
$(\pi, {\al H.})$ of the \mbox{*--algebra} ${\al R.}_0$  
generated by the resolvents $R(\lambda,f)$ with 
$\lambda \in \R \backslash \{0\}$, $f \in X$, 
putting in analogy to relation~(\ref{RepEq})
\begin{equation} \label{RepDef}
\pi(R(\lambda,f)) \doteq
\begin{cases}
\ \pi_0 \, (R(\lambda + i\varphi(f_Z) \, ,f_N)) & \mbox{if}
\quad  f = f_Z + f_N \in Y  \\ 
\ 0 & \mbox{if} \quad f \in X \backslash Y 
\end{cases}
\end{equation} 
and extending this definition of $\pi$ to arbitrary finite 
sums and products of resolvents by linearity and multiplicativity. 
This definition is consistent with relations
(\ref{Resolv}) to (\ref{Riden}) 
and thus determines a representation of ${\al R.}_0$
on $\al H.$ which can be extended to $\rsl$ by continuity. 
By construction, $\pi(\rsl) = \pi_0({\al R.}(Q^\perp \cap
Y, \sigma))$, hence $(\pi,{\al H.})$ is also irreducible. 

It remains to show that the element $(X_R, \chi_T) \in I (X,\sigma)$ 
associated with $(\pi,{\al H.})$ coincides with the given $(Y,\chi)$. 
Since $(\pi_0,{\al H.})$ is a regular representation of 
${\al R.}(Q^\perp \cap Y, \sigma)$, it follows from the defining
relation (\ref{RepDef}) that $\ker\big(\pi(R(\lambda,f))\big) = \{0\}$ 
iff $f \in Y$, hence $X_R = Y$. Moreover, by Lemma 
\ref{Xdecomp}(iii) $X_T = X_R \cap {X_R}^\perp = Y \cap Y^\perp = Z$
and there exists a pure state $\chi_T$ on the abelian algebra 
$\al R._T = {C^*\{\rlf :  \lambda \in \Ro  \, , \, f\in X_T = Z \}} = \al Z.$
such that $\pi \rest {\al R._T} = \chi_T(\, \cdot \,) \, \un_{\al H.}$.
Now according to relation (\ref{RepDef})
\[ \pi(R(\lambda,f_Z)) = (i \lambda - \varphi(f_Z))^{-1} \, 
\un_{\al H.} = \chi(R(\lambda,f_Z)) \, \un_{\al H.} \, , 
\quad f_Z \in Z = X_T \, ,
\]
and $\chi$ is, by assumption, a pure state on 
$\al Z. = \al R._T$. Hence $\chi = \chi_T$, completing the proof.
\end{beweis}

Note that if $(X,\,\sigma)$ is infinite dimensional, 
the map $\iota$ is definitely not injective, as there exist 
inequivalent regular irreducible representations. There is then  
no simple characterization of the spectrum of $\rsl$.  
Focusing on the finite dimensional case, we determine
next the primitive ideals of $\rsl$, \ie the kernels of its 
irreducible representations. Before entering into this analysis
let us recall two basic facts about (throughout this paper 
always closed and two--sided)
ideals, cf.\ \cite[Add.~1.9.12]{Dix} and \cite[Prop.~II.8.2.4]{Bla}.

\begin{lem}
\label{InductIdeal}
Let $\al A.$ be a C*-algebra,
\begin{itemize}
\item[(i)]
Given ideals $\al J._1$ and $\al J._2$ of
$\al A.$, then $\al J._1\cap\al J._2=\al J._1\cdot\al J._2=
\al J._2\cdot\al J._1$.
\item[(ii)]
Let $\{\al A._i : i\in I\}$ be a set of C*-subalgebras
such that $\al A._0\doteq \bigcup\limits_{i\in I}\al A._i$ 
is a dense *-subalgebra of $\al A.$.
If $\al J.$ is any
ideal of $\al A.$, then 
$\al J.\cap\al A._0 = \bigcup\limits_{i\in I}\big(\al A._i \cap
\al J.\big)$ is dense in $\al J.$.
\end{itemize}
\end{lem}

We can now prove:

\begin{pro}
\label{PrimIdeals}
Let $(X,\,\sigma)$ be a  finite dimensional symplectic space and
let $\widehat{\pi} \in \widehat{\rsl}$ 
with associated pair $\iota(\widehat{\pi})={(X_R,\chi_T)}\in
I (X,\sigma)$. Then
\begin{itemize}
\item[(i)]
$\ker \widehat{\pi}$ (\ie the common kernel of all 
irreducible representations $(\pi, \al H.)$ 
belonging to the class $\widehat{\pi}$)
is the (possibly zero) ideal generated in $\rsl$ by
\[
S(X_R,\chi_T) \doteq
\{R(\lambda,f) : f\in X_S \} \, {\textstyle \bigcup} \, 
\{(R(\mu,g)-\chi_T(R(\mu,g)) \, \un) : g \in X_T \} \, . 
\] 
This ideal does not depend on the choice of $\lambda, \mu \in \Ro$. 
\item[(ii)]
The map $\widehat{\pi} \mapsto \ker \widehat{\pi}$ is a bijection 
from $\widehat{\rsl}$ to the set of primitive ideals  of 
$\rsl$.
\end{itemize} 
\end{pro}

\begin{beweis}
(i) Let $\al J.(X_R,\chi_T) \subset \rsl$ denote the 
ideal generated by the set $S(X_R,\chi_T)$. It follows
by a routine computation from
the resolvent equation (\ref{Resolv}) and the fact that 
$\chi_T$ is a pure state on the abelian algebra 
$\al R._T = C^*\{R(\mu,g) : \mu \in \Ro \, , \, g \in X_T \}$
that $\al J.(X_R,\chi_T)$ does not depend 
on the  choice of $\lambda, \mu \in \Ro$. 

Let $(\pi, \al H.)$ be any irreducible representation 
in the class of the given $\widehat{\pi}$. By the very 
definition of the pair $(X_R, \chi_T)$ the set 
$S(X_R,\chi_T)$ lies in the kernel of $\pi$
for any \mbox{$\lambda, \mu \in \Ro$},
hence $\al J.(X_R,\chi_T) \subseteq \ker \pi$. 
For the proof that one has 
equality, we choose a space $Q$ as in relation~(\ref{Qdecomp}) 
and make use of the unique decomposition of 
$X_R$ into $X_T = X_R  \cap {X_R}^\perp$ and $X_N = Q^\perp \cap X_R$. Let
$f_R = f_T + f_N \in X_R$ with $f_T \in X_T$, $f_N \in X_N$
and let $\lambda, \mu \in \Ro$, $\lambda \neq \mu$. It follows
from the resolvent equation (\ref{Resolv}) that 
$(\un - i \mu \, R(\lambda,f_T))^{-1} = (\un + i \mu \, R(\lambda - \mu,f_T))$
and taking also into account that $\sigma(f_R,f_T) = 0$ 
equation (\ref{Rsum}) implies
\[
R(\lambda, f_R) - R(\mu, f_N)  
= \big( \un - i\mu \, R(\lambda,f_T) \big) \,
\big(R(\lambda,f_R) - R(\mu, f_N) 
- i\mu \, R(\lambda, f_R) R(\mu,f_N) \big)   \, .
\]
Because of the analyticity properties of the 
resolvents this equality extends to  
$\mu \in \C \backslash i \, \R$, hence putting 
$\mu = -i \, \chi_T(R(\lambda,f_T))^{-1}$ one obtains 
$( R(\lambda, f_R) - R(\mu, f_N) ) \in \al J.(X_R,\chi_T)$. 
Since also \mbox{$R(\lambda, f_S) \in \al J.(X_R,\chi_T)$} for 
$f_s \in X_S$ each element of the dense 
subalgebra ${\al R.}_0 \subset \rsl$ generated
by all polynomials in $\rlf$ for $\lambda \in \Ro$, $f \in X$
can be decomposed into a sum of elements of 
$\al R.(X_N, \sigma)$ and of $\al J.(X_R,\chi_T)$. 
Now the 
representation $\pi \rest  \al R.(X_N,\sigma)$ is, by the very 
definition of the symplectic subspace $X_N$, regular and 
hence faithful according to Proposition \ref{UniqueR}(iii) 
and $\pi \rest \al J.(X_R,\chi_T) = 0$. So
this decomposition  extends by the continuity of $\pi$ 
uniquely to all elements
of $\rsl$. Consequently $\ker {\pi} \subseteq  \al J.(X_R,\chi_T)$, 
completing the proof of part (i).  

(ii) By definition, the map $\widehat{\pi} \to \ker \widehat{\pi}$ 
is a surjection from $\widehat{\rsl}$ to the set of primitive 
ideals of $\rsl$, so we only need to prove injectivity. 
Let $\widehat{\pi},\,\widehat{\pi}^\prime \in \widehat{\rsl}$ 
such that  
\mbox{$\ker\widehat{\pi} = \ker\widehat{\pi}^\prime$} and let
$\iota(\widehat{\pi})={(X_R,\chi_T)}$,  \ 
$\iota(\widehat{\pi}^\prime) = 
(X_R^\prime ,\chi_T^{\, \prime}) $ 
be the associated pairs in $I (X,\sigma)$.
We pick representations 
$(\pi, \al H.)$,  $({\pi}^\prime, {\al H.}^\prime)$ in the classes of 
of $\widehat{\pi}$ and $\widehat{\pi}^\prime$, respectively.
According to Lemma \ref{Xdecomp}
$\pi\big(\rlf\big) \neq \{0\}$ iff $f\in X_R$ and, similarly, 
$\pi^\prime\big(\rlf\big) \neq \{0\} $  iff $f\in {X_R}^\prime$. 
Since $\ker \pi = \ker \pi^\prime$ it follows that 
$X_R =  {X_R}^\prime$, hence $X_T ={X_T}^\prime$. Moreover, 
in view of the inclusion  
$S(X_R, \chi_T) \subset \ker \pi = \ker \pi^\prime$ we have 
\[
\big(\chi_T^{\, \prime}(R(\lambda,g)) - \chi_T(R(\lambda,g)) \big) 
\, \un_{{\al H.}^\prime} =
\pi^\prime (R(\lambda,g) - \chi_T(R(\lambda,g)) \, \un ) = 0 \, ,
\quad g \in X_T \, .
\] 
\ie $\chi_T = \chi_T^{\, \prime}$ and 
consequently  $I (\widehat{\pi}) =
I (\widehat{\pi}^\prime)$. According to 
Proposition \ref{PrimIdeals}(i) this implies 
$\widehat{\pi} = \widehat{\pi}^\prime$, 
completing the proof. 
\end{beweis} 

Property (ii) in this proposition 
is a remarkable feature of the resolvent algebra,
shared with the abelian C*-algebras. It rarely holds for 
noncommutative algebras. Consider for example the Weyl algebra 
which, being simple, has only faithful representations. 
As a matter of fact, property (ii) does not hold either for the resolvent 
algebra if the underlying symplectic space is infinite dimensional,
cf.\ the remark made after Proposition \ref{Irrep}. 

Next, we consider the partial order of the set of primitive 
ideals in $\rsl$, 
$\ker \widehat{\pi} \subseteq \ker \widehat{\pi}^\prime$, where 
strict inclusions will be denoted by $\ker \widehat{\pi} \subsetneq 
\ker \widehat{\pi}^\prime$. 
\begin{defi}
Let $\al A.$ be a unital C*--algebra and let $\widehat{\al A.}$ be the
corresponding set of equivalence classes of irreducible
representations of $\al A.$. The maximal length of 
strictly increasing chains of 
(possibly zero) primitive ideals in ${\al A.}$ 
is denoted by $\, L({\al A.})$, \ie
\[
L({\al A.}) \doteq 
\sup \, \{ n \in \N : \ker \widehat{\pi}_1 \subsetneq 
\dots  \subsetneq \ker \widehat{\pi}_n \, , \ 
\widehat{\pi}_1, \dots , \widehat{\pi}_n  \in \widehat{\al A.} \} 
\]
and $L({\al A.}) \doteq \infty$ in case the supremum does 
not exist. The quantity $L({\al A.})$ is clearly an 
isomorphism invariant of C*--algebras. 
\end{defi}

We are now in a position to prove the main result of this section.
\begin{teo}
\label{DimThm}
Let $(X,\,\sigma)$ be a symplectic space of arbitrary dimension. Then
\begin{itemize}
\item[(i)]  $L(\rsl)= 
\dim(X)/2  +1 \,$ if $\, \dim(X) < \infty \,$ 
and $\, L(\rsl) = \infty \,$ otherwise.  
\item[(ii)] The isomorphism classes of the set of all 
resolvent algebras associated with finite dimensional symplectic 
spaces are completely characterized by $L(\rsl)$.
\end{itemize}
\end{teo}
\begin{beweis} (i) First, let $X$ be finite dimensional, let 
$(\pi, \al H.)$, $(\pi^\prime, {\al H.}^\prime)$ be irreducible 
representations in the classes 
$\widehat{\pi}, \widehat{\pi}^\prime \in \widehat{\rsl}$,
respectively, 
and let $\ker \widehat{\pi} \subseteq \ker \widehat{\pi}^\prime$. 
According to Proposition~\ref{PrimIdeals}(i) 
the kernels of the two representations coincide with the 
two--sided ideals generated by the corresponding 
sets $S(X_R,\chi_T)$ and $S({X_R}^\prime, {\chi_T}^\prime)$,
respectively. Now if $f \in X_S = X \backslash X_R$, 
\ie $\pi(R(\lambda,f)) = 0$ for $\lambda \in \Ro$, one 
also has $\pi^\prime(R(\lambda,f)) = 0$, 
\ie $f \in {X_S}^\prime = X \backslash {X_R}^\prime$,
proving $X_S \subseteq {X_S}^\prime$
which implies $X_R \supseteq {X_R}^\prime$. Similarly, if
$g \in X_T  = X_R \cap {X_R}^\perp $ one has 
$\pi\big(R(\mu,g) - \chi_T(R(\mu,g)) \, \un \big)
= 0$ for $\mu \in \Ro$, hence 
$\pi^\prime \big(R(\mu,g) - \chi_T(R(\mu,g)) \, \un \big) = 0$. 
But $\chi_T(R(\mu,g)) \neq 0$, so it follows that
$g \in {X_T}^\prime$, hence $X_T \subseteq {X_T}^\prime$, 
 and ${\chi_T}^\prime \rest  {X_T} = \chi_T$. Moreover, in case of 
a strict inclusion $\ker \widehat{\pi} \subsetneq \ker \widehat{\pi}^\prime$
one has $X_R \supsetneq {X_R}^\prime$; for otherwise 
$S(X_R,\chi_T) = S({X_R}^\prime, {\chi_T}^\prime)$ in
conflict with Proposition \ref{PrimIdeals}. 

Now according to Lemma~\ref{Xdecomp}(iv)  
one can always extend the space $X_T = X_R \cap {X_R}^\perp $ to a 
nondegenerate subspace $Q \subset X$ by augmenting it with some 
complementary space $X_T^{\, \sim}$ of the same dimension. The space 
$X_N \doteq Q^\perp \cap X_R = X_T^{\, \sim \, \perp} \cap X_R$ is 
nondegenerate by construction and has the dimension 
$\dim(X_N) = \big(\dim(X_R) - \dim(X_T)\big)$ which is even 
and bounded by $0 \leq \dim(X_N) \leq \dim(X)$. 
The same statements apply {\it mutatis mutandis} to the 
spaces ${X_R}^\prime$,~${X_T}^\prime$ and ${X_N}^\prime$ 
affiliated with the second representation. In case of a 
strict inclusion 
\mbox{$\ker \widehat{\pi} \subsetneq \ker \widehat{\pi}^\prime$}  
it follows from the preceding discussion that 
$X_R \supsetneq {X_R}^\prime$. Since 
$X_T \subseteq {X_T}^\prime$ this implies $\dim(X_N) > \dim({X_N}^\prime)$,
so the length $l$ of
any given strictly increasing chain of primitive ideals in 
$\rsl$ complies with the upper bound $l \leq \dim(X)/2 + 1$. 
In order to exhibit a chain where one has equality we 
choose any strictly decreasing sequence
of nondegenerate subspaces $X_0 = X \supsetneq X_1 \dots 
\supsetneq X_{\dim(X)/2} = \{0\}$ and consider the family of 
representations $(\pi_n, {\al H.}_n)$ 
of $\rsl$ for which $\pi_n \rest {\al R.} (X_n, \sigma)$
is regular, acts irreducibly on ${\al H.}_n$ and 
$\pi_n(R(\lambda,f)) = 0$, $f  \in X \backslash X_n$ for 
$n = 0, \dots,  \dim(X)/2$, cf.\ definition~(\ref{RepDef}).
The resulting strictly increasing chain of primitive ideals
$\ker \pi_n$, $n = 0, \dots,  \dim(X)/2$, has the 
desired length, 
hence $L(\rsl) = \dim(X)/2 + 1$, as claimed. It also
follows from the latter construction that $L(\rsl) = \infty$
if $X$ is infinite dimensional, thereby completing the proof of the first
part of the statement.

(ii) As we have seen, $L$ defines an isomorphism 
invariant from which the dimension of the symplectic 
space underlying a resolvent algebra can be recovered. 
Conversely, let $(X,\sigma)$, $(X^\prime, \sigma^\prime)$ be 
symplectic spaces of equal finite dimension and let 
$\rsl$, ${\al R.} (X^\prime, \sigma^\prime)$ be the 
corresponding resolvent algebras. There then exists
a symplectic transformation $\Gamma : X \rightarrow X^\prime$
mapping the first space onto the second one and satisfying
$\sigma^\prime(\Gamma f, \Gamma g) = \sigma(f,g)$, \ $f,g \in X$. 
This transformation induces a bijection between the generating 
elements of the resolvent algebras given by 
$\gamma(R(\lambda,f)) \doteq R^{\, \prime}(\lambda,\Gamma f)$, cf.
\cite[Thm.~5.3(ii)]{BuGr2}, 
which is compatible with the defining relations. It therefore
extends to an isomorphism $\gamma : \rsl \rightarrow {\al R.}
(X^\prime, \sigma^\prime)$, as claimed. Thus $L$ 
provides a complete algebraic invariant for the family of 
resolvent algebras associated with finite dimensional 
symplectic spaces.
\end{beweis}

Another consequence of the preceding results is the following 
theorem which throws further light on the structure of 
resolvent algebras. It ought to be mentioned in this context that
the resolvent algebras are not separable \cite[Thm.~5.3]{BuGr2}.

\begin{teo}
\label{Rnuclear}
Let $(X,\sigma)$ be a symplectic space of arbitrary dimension. Then
\begin{itemize}
\item[(i)]
 $\rsl$ is a nuclear C*-algebra,
\item[(ii)] $\rsl$ is a Type~I (postliminal) C*-algebra
iff  ${\rm dim}(X)<\infty$.
\end{itemize}
\end{teo}
\begin{beweis}
As Type~I C*-algebras are nuclear, it follows from part~(ii) and 
Proposition~\ref{UniqueR}(ii)
that $\rsl$ is an inductive limit of nuclear algebras, hence nuclear,
cf.\ \cite[Prop.~IV.3.1.9]{Bla}. Thus we only need to prove part~(ii).

By  Theorem~IV.1.5.7 and  IV.1.5.8 (nonseparable case) 
in~\cite{Bla} we know that $\rsl$ is a Type~I C*-algebra iff it is
  GCR, \ie its image for each irreducible representation has 
nonzero intersection with the compacts \cite[Def.~IV.1.3.1]{Bla}.
Let ${\rm dim}(X)<\infty$, then  for an irreducible representation 
$(\pi, \al H.)$ we know from relation (\ref{RepEq}) 
and the subsequent remarks that
$\pi\big(\rsl\big)=\pi\big(\al R.(X_N,\sigma)$,
where $X_N \subseteq X$ is a nondegenerate subspace. Moreover, 
\mbox{$\pi \rest \al R.(X_N,\sigma)$} is a regular irreducible 
representation, hence its range contains the 
compacts according to \cite[Thm.~5.4(i)]{BuGr2}. 
So  $\rsl$  is GCR  hence Type I if $X$ is finite dimensional. 
Conversely, if $(\pi, \al H.)$ is a faithful irreducible representation
of $\rsl$ such that $\pi(\rsl)$ contains the compacts, then their respective
pre--image constitutes a non--zero ideal in $\rsl$
which is minimal in view of Lemma \ref{InductIdeal}(i). 
But as we will show in 
Theorem~\ref{MinId}(ii) there are no such ideals if $\dim(X) = \infty$,
hence $\rsl$ is not GCR in this case.
\end{beweis}

\section{Principal ideals}
\setcounter{equation}{0}
\label{PrincipalIdeals}

In the preceding section we have characterized all the primitive ideals of
$\rsl$ for the case $\dim(X) < \infty$. If $\dim(X) = \infty$, 
an exhaustive characterization of these ideals seems
a hopeless task, however. The concept of principal ideals, 
\ie ideals generated by some element of the algebra, is more 
appropriate then for structural analysis. 

Throughout this section, we consider 
symplectic spaces $(X,\sigma)$ of arbitrary dimension
unless otherwise stated. 
Let $\al J. \subset \rsl$ be any ideal. Since $\rsl$ is the
C*--inductive limit of the algebras $\al R.(S,\sigma)$ based on all 
finite dimensional nondegenerate subspaces $S \subset X$, it 
follows from Lemma~\ref{InductIdeal}(ii) that  $\al J.$ is the
C*--inductive limit of the ideals 
\mbox{$\al J.(S) \doteq \al J. \cap
\al R.(S,\sigma)$} in $\al R.(S,\sigma)$. As we have seen, the 
latter ideals are built from principal ideals generated by
the operators $(R(\lambda, f) - \rho \, \un)$, 
where $\rho$ belongs to $\spec (R(\lambda, f))$, the spectrum of the 
operator $\rlf$ which 
according to relations (\ref{Rinvol}) and (\ref{Rccr}) 
is normal. Hence these principal ideals are building blocks of all ideals.
It is therefore warranted to have a closer look at their structure.
We begin with a preparatory lemma which slightly generalizes 
Theorem~4.1(iv) and Proposition~8.1(ii) in \cite{BuGr2}. 
\begin{lem} 
\label{SingularStates}
Let $f \in X \backslash \{0\}$, $\lambda \in \Ro$, 
and let $\rho   \in  \spec(R(\lambda, f))$. 
\begin{itemize}
\item[(i)] There exists a pure state $\omega$ on $\rsl$ such that
$(R(\lambda, f) - \rho \,  \un) \in \ker \, \omega$. 
\item[(ii)] Given a state $\omega$ on $\rsl$ such that
$(R(\lambda, f) - \rho \,  \un) \in \ker \, \omega$,
then \mbox{$\pi_\omega(R(\lambda, f) - \rho \,  \un) =0$}, where 
$(\pi_\omega, {\al H.}_\omega)$ denotes the GNS representation
induced by $\omega$. Moreover, if $\rho \neq 0$ then
$\pi_\omega(R(\mu,g)) = 0$ for all $\mu \in \R \backslash \{0\}$ and
$g \in X$ satisfying $\sigma(g,f) \neq 0$.
\end{itemize} 
\end{lem}

\begin{beweis}
(i) Since $\rho$ is contained in the spectrum of 
$R(\lambda,f)$, there exists a pure state 
$\chi$ on the abelian C*--algebra generated by 
$R(\mu,f)$, $\mu \in \Ro$ and $\un$
such that $\chi(R(\lambda,f) - \rho \, \un)=0$.
By the Hahn--Banach theorem, $\chi$ can be
extended to a pure state $\omega$ on $\rsl$.  \\
(ii) Since $R(\lambda,f)$ is a normal operator 
it follows from the resolvent equations 
(\ref{Resolv}), (\ref{Rinvol})  that its spectral values $\rho$
satisfy 
$\overline{\rho} - \rho = 2 i \lambda \, |\rho|^2$. Hence one obtains for
the given $\omega$
\begin{equation*}
\omega\big((R(\lambda,f) - \rho \un)^*(R(\lambda,f) - \rho \un))
= (2i \lambda)^{-1} \, \omega(R(\lambda,f)^* - R(\lambda,f))
- |\rho|^2 = 0 \, .
\end{equation*}
Let $\Omega_\omega \in {\al H.}_\omega$ be the 
GNS--vector derived from $\omega$. The preceding equality  
and the fact that $R(\lambda,f)$ is normal implies 
$\pi_\omega(R(\lambda,f) - \rho \, \un) \, \Omega_\omega = 0 
= \pi_\omega(R(\lambda,f) - {\rho} \, \un)^* \, \Omega_\omega $. 
It will be shown in two steps that these equalities imply 
$\ker \pi_\omega (R(\lambda,f) - \rho \, \un) = {\al
  H.}_\omega$. First, let $\rho = 0$. It 
then follows from relation~(\ref{Rccr}) that 
$\ker \pi_\omega (R(\lambda,f))$ is stable under the action of 
$\pi_\omega(R(\mu,g))$ for any $\mu \in \Ro$, $g \in X$.
Since $\Omega_\omega \in \ker \pi_\omega (R(\lambda,f))$ is cyclic 
for $\pi_\omega(\rsl)$ this implies $\pi_\omega(R(\lambda, f)) = 0$. Second, if
$\rho \neq 0$ it is still true that 
$\ker \pi_\omega (R(\lambda,f) - \rho \, \un)$
is stable under the action of the operators 
$\pi_\omega(R(\mu,g))$ whenever $\sigma(f,g) = 0$. So let $g \in X$ 
be such that $\sigma(f,g) \neq 0$ and let 
$\mu \in \Ro$. By relation~(\ref{Rccr})  
\begin{equation*}
0 = \omega([R(\lambda,f), R(\mu,g)]) =
\sigma(f,g) \, \rho^2 \omega(R(\mu,g)^2) \, ,
\end{equation*}
where it has been used that 
$\omega(R(\lambda,f) A) = \omega(A R(\lambda,f) ) = \rho \, \omega(A)$
for $A \in \rsl$. Hence $\omega(R(\mu,g)^2)= 0$.
Since 
$\mu \mapsto R(\mu,g)$ is differentiable and 
$ i {d \over d \mu} \, R(\mu,g) = R(\mu,g)^2 $
according to relation (\ref{Resolv}) one arrives at 
$i {d \over d \mu} \, \omega(R(\mu,g)) = 0$. Thus 
$\mu \mapsto \omega(R(\mu,g)) = \mbox{\it const} = 0$, 
where the second equality follows from the 
bound $\| R(\mu,g) \| \leq |\mu|^{-1}$ for 
large $\mu$. According to 
the first step this entails $\pi_\omega(R(\mu,g)) = 0$, hence  
$\ker \pi_\omega (R(\lambda,f) - \rho \, \un)$ is stable under the 
action of all operators $\pi_\omega(R(\mu,g))$ with 
$\mu \in \R \backslash \{0\}$, $g \in X$. 
Since $\Omega_\omega \in  \ker \pi_\omega (R(\lambda,f) - \rho \,
\un)$ it is then clear that 
$\pi_\omega (R(\lambda,f) - \rho \, \un) = 0$ for $\rho \neq 0$
as well. The last part of the statement has already been 
established in the preceding step, completing the proof. 
\end{beweis}

\begin{pro} \label{PrincId}
Let $\lambda \in \R \backslash \{0\}$, 
$f \in X \backslash \{0\}$, $\rho \in  \spec(R(\lambda, f))$  and let 
${\al I.}(\lambda,f;\rho)$ be the ideal generated by 
$(R(\lambda, f) - \rho \,  \un)$, \viz  
\mbox{${\al I.}(\lambda,f;\rho) 
\doteq [\rsl (R(\lambda, f) - \rho \,  \un) \rsl]$}.\footnote{Here
  and in the following $[ \ \cdot \ ]$ denotes the closed linear span of 
its argument}
\begin{itemize}
\item[(i)] $\al I.(\lambda,f;\rho)$ is proper
\item[(ii)]  $\al I.(\lambda,f;\rho) = [\rsl (R(\lambda, f) - \rho \,  \un) ]
=  [(R(\lambda, f) - \rho \,  \un) \rsl]$   
\item[(iii)] $(R(\mu,f) - 
\mbox{\large$\frac{\rho}{1+i \, (\mu - \lambda) \, \rho}$} \, \un) 
\in \al I.(\lambda,f;\rho)$ 
for all $\mu \in \R \backslash \{0\}$  
\item[(iv)] If $\rho \neq 0$ then
$R(\mu,g) \in \al I.(\lambda,f;\rho)$ for all $\mu \in \R \backslash \{0\}$ and
$g \in X$ such that $\sigma(g,f) \neq 0$.
\end{itemize}
\end{pro}

\begin{beweis}
(i) According to part (i) of the preceding lemma there exists 
a state $\omega$ on $\rsl$ such that 
$(R(\lambda,f) - \rho \, \un) \in \ker \omega$. 
Hence, by part (ii) of the lemma, 
\mbox{$\pi_\omega(R(\lambda,f) - \rho \, \un) = 0$}.
Thus ${\al I.}(\lambda,f;\rho)$, being non--trivial and 
contained in the kernel of the representation 
$\pi_\omega$, is a proper ideal. 

(ii) Any closed left ideal is the intersection of the left kernels
of all states which contain it, cf.~\mbox{\cite[Lem.~3.13.5]{Ped}}. 
Let $\omega$ be any state on $\rsl$ such that
\mbox{$[\rsl (R(\lambda,f) - \rho \, \un)] \subseteq {\al L.}_{\, \omega}$}, 
where ${\al L.}_{\, \omega}$ denotes the left kernel of $\omega$.
Then 
$[\rsl (R(\lambda,f) - \rho \, \un) \rsl ] \subseteq {\al L.}_{\,
  \omega}$ according to Lemma \ref{SingularStates}(ii). 
Taking intersections with 
regard to all such states $\omega$, one arrives 
at~$[\rsl (R(\lambda,f) - \rho \, \un) \rsl ] \subseteq 
[\rsl (R(\lambda,f) - \rho \, \un)]$. The opposite inclusion is
trivial, hence 
$[\rsl (R(\lambda,f) - \rho \, \un)] =
[\rsl (R(\lambda,f) - \rho \, \un) \rsl ] $.
Replacing in this equality $\rlf$ by $\rlf^* = R(-\lambda,f)$, \ 
$\rho$ by its complex conjugate $\overline{\rho}$ and recalling 
that any closed two--sided ideal of a 
unital C*--algebra is stable under taking adjoints 
one also obtains 
$[(R(\lambda,f) - \rho \, \un) \rsl  ] =
[\rsl (R(\lambda,f) - \rho \, \un) \rsl ] $, as claimed. 

(iii) Relation (\ref{Resolv}) implies 
$(1 + i(\mu-\lambda) \, \rho) \, R(\mu,f) - \rho \, \un 
= (\un + i(\lambda - \mu) \, R(\mu,f)) (R(\lambda,f) - \rho \, \un)$
from which the assertion follows. 

(iv) If $\rho \neq 0$ and $\omega$ is any state such that
$[\rsl (R(\lambda,f) - \rho \, \un)] \subseteq {\al L.}_{\, \omega}$
it follows from the last part of Lemma \ref{SingularStates}(ii) that
$R(\mu,g) \in {\al L.}_{\, \omega}$ for $\mu \in \R \backslash \{0\}$
and any $g \in X$ satisfying $\sigma(f,g) \neq 0$. Consequently
$R(\mu,g) \in [\rsl (R(\lambda,f) - \rho \, \un)] $. 
\end{beweis}

It is apparent from the preceding discussion that the 
principal ideals under consideration are in one--to--one correspondence with 
representations in which the underlying generators $\phi(f)$
of the resolvents 
have sharp values (including the singular case $\infty$).
We turn now to the analysis of intersections of these
principal ideals which turn out to be principal ideals as well.
Before proving this we need to establish  
the following lemma.

\begin{lem}
Let $f \in X \backslash \{0\}$, $\lambda \in \R \backslash \{0\}$
and $\rho \in  \spec(R(\lambda, f))$,  and let 
${\al I.}(\lambda,f;\rho)$ be the ideal defined in the preceding
proposition. For any ideal $\al J. \subset \rsl$ one has 
\begin{equation*}
{\al I.}(\lambda,f;\rho) \, {\textstyle \bigcap} \, \al J. =
[(R(\lambda,f) - \rho \, \un) \al J.] =
[\al J. (R(\lambda,f) - \rho \, \un)] =
[\al J. (R(\lambda,f) - \rho \, \un) \al J.] \, .
\end{equation*}
\end{lem}
\begin{beweis}
According to Lemma \ref{InductIdeal}(i) \  
 ${\al I.}(\lambda,f;\rho) \cap \al J. = 
{\al I.}(\lambda,f;\rho) \cdot \al J. $. Hence, making use of 
part~(ii) of the preceding proposition, one obtains 
\begin{equation*}
{\al I.}(\lambda,f;\rho) \, {\textstyle \bigcap} \, \al J. =
[(R(\lambda,f) - \rho \, \un) \rsl] \cdot {\al J.}
\subseteq [(R(\lambda,f) - \rho \, \un) \rsl {\al J.}]
= [(R(\lambda,f) - \rho \, \un) {\al J.}] \, .
\end{equation*}
Clearly $[(R(\lambda,f) - \rho \, \un) {\al J.}] \subseteq 
{\al I.}(\lambda,f;\rho) \, {\textstyle \bigcap} \, \al J.$, 
so the first equality in the statement follows and in a 
similar manner one obtains the second equality. For the last 
equality one makes use of the fact that 
$[\al J. (R(\lambda,f) - \rho \, \un)] =
{\al I.}(\lambda,f;\rho) \, {\textstyle \bigcap} \, \al J. $
is an ideal, hence 
\begin{equation*}
[\al J. (R(\lambda,f) - \rho \, \un) \al J.] 
\subseteq [[\al J. (R(\lambda,f) - \rho \, \un) ] \, \al J. ] 
\subseteq  [\al J. (R(\lambda,f) - \rho \, \un) ] \, ,
\end{equation*}
and the opposite inclusion holds since 
any ideal $\al J.$ has an approximate identity.
\end{beweis}

\begin{pro} \label{IdIntersect}
Let 
$f_m \in X \backslash \{0\}$, $\lambda_m \in \R \backslash \{0\}$,  
$\rho_m \in  \spec(R(\lambda_m, f_m))$ and let 
${\al I.}(\lambda_m,f_m;\rho_m)$ be the corresponding ideals introduced
in Proposition \ref{PrincId}, $m = 1, \dots , n$. Then
\begin{equation*}
\begin{split}
& {\textstyle \bigcap_{m=1}^n} \, {\al I.}(\lambda_m,f_m;\rho_m) \ = \ 
[\rsl \ \Pi_{m=1}^n \big( R(\lambda_m, f_m) - \rho_m \, \un \big) \ \rsl]
\\
& = [\, \Pi_{m=1}^n \big( R(\lambda_m, f_m) - \rho_m \, \un \big) \ \rsl]
= [\rsl \ \Pi_{m=1}^n \big( R(\lambda_m, f_m) - \rho_m \, \un \big) ]
\, , 
\end{split}
\end{equation*}
where the order of the operators $\big( R(\lambda_m, f_m) - \rho_m \, \un
\big)$ in the products is arbitrary. 
\end{pro}

\pagebreak
\begin{beweis}
The proof proceeds by induction in $n$. For $n=1$ the statement   
was established in Proposition \ref{PrincId}. Putting
${\al J.}_n \doteq  {\textstyle \bigcap_{m=1}^n} \, 
{\al I.}(\lambda_m,f_m;\rho_m) $ it follows from the preceding lemma
and the induction hypothesis that 
\begin{equation*}
\begin{split}
{\al J.}_{n+1} & = 
{\al I.}(\lambda_{n+1},f_{n+1};\rho_{n+1}) \, {\textstyle \bigcap} \,
{\al J.}_n  
  = [\big( R(\lambda_{n+1}, f_{n+1}) - \rho_{n+1} \, \un \big) \ {\al
   J.}_n  ] \\
& = [\big( R(\lambda_{n+1}, f_{n+1}) - \rho_{n+1} \, \un \big) \,
\Pi_{m=1}^{n} \big( R(\lambda_m, f_m) - \rho_m \, \un \big) \, \rsl] \, .
\end{split}
\end{equation*}
Similarly 
${\al J.}_{n+1} =  [\rsl \ \Pi_{m=1}^{n} \big( R(\lambda_m, f_m) - \rho_m
\, \un \big) \, \big( R(\lambda_{n+1}, f_{n+1}) - \rho_{n+1} \, \un
\big)]$. 
Moreover, since ${\al J.}_{n+1}$ is an ideal, the preceding 
equality implies 
\[ {\al J.}_{n+1} =  [\rsl \ \Pi_{m=1}^{n} \big( R(\lambda_m, f_m) - \rho_m
\, \un \big) \, \big( R(\lambda_{n+1}, f_{n+1}) - \rho_{n+1} \, \un
\big) \rsl] \, .
\]
Since the intersection of sets is stable under their permutation 
the order of the factors $\big( R(\lambda_\cdot, f_\cdot) - \rho_\cdot \, \un
\big)$ in the products is arbitrary, completing the proof of 
the statement. 
\end{beweis}

As already mentioned at the end of Sec.~\ref{ResBasics}  
the intersection of any two nonzero ideals in $\rsl$
is nonzero. The situation changes, however, if one proceeds to infinite 
intersections. There one encounters a marked difference between 
the resolvent algebras of finite and of infinite systems. 

\medskip
\begin{teo} \label{MinId}
Let $\al J. \subset \rsl$ be the intersection of all nonzero ideals
of $\rsl$.
\begin{itemize}
\item[(i)] If $\dim(X) < \infty$ then the ideal $\al J.$ is isomorphic to the 
C*--algebra $\al K.$ of compact operators.  
\item[(ii)] If  $\dim(X) = \infty$ then $\al J. = \{0\}$. In fact,
there exists no nonzero minimal ideal of $\rsl$ in this case.  
\end{itemize}
\end{teo}
\begin{beweis}
(i) According to Theorem \ref{Rnuclear}(ii) the algebra $\rsl$ is of 
type I (postliminal) if $\dim(X) < \infty$. Picking 
any faithful irreducible  representation  $(\pi, \al H.)$ 
of  $\rsl$ one has  
$\al K. \subset \pi(\rsl)$, so the pre--image 
$\pi^{-1}(\al K.)$ is a nonzero ideal of $\rsl$. Given any 
other nonzero ideal $\al I. \subset \rsl$ one has 
$\al I. \cap \pi^{-1}(\al K.) \neq \{0\}$ since $\rsl$ is prime
and consequently 
$\pi(\al I.) \cap \al K. \neq \{0\}$ since $\pi$
is faithful. So the ideal $\pi(\al I.) \subset \al B.(\al H.)$ contains some
nontrivial compact operator and consequently 
$\al K. \subseteq \pi(\al I.)$ since $\pi$ is irreducible. 
Hence $\pi^{-1}(\al K.) \subseteq \al I.$ for any non--zero ideal $\al
I.$, proving the first statement. 

(ii) Let $\al J.$ be a nonzero minimal ideal of $\rsl$. It follows
from Proposition~\ref{UniqueR}(ii)  and 
Lemma~\ref{InductIdeal}(ii) that there exists a 
finite dimensional non--degenerate subspace $S_0 \subset X$ such that 
$\al J. \cap \al R. (S_0,\sigma) \neq \{0\}$. Hence 
$\al J. \cap \al R. (S_0,\sigma)$ is a nonzero 
ideal of $\al R. (S_0,\sigma)$ which, according to the first
step, contains a distinguished algebra
$\al K. (S_0)$ which is isomorphic to the algebra of compact
operators. Hence $[\rsl \al K. (S_0) \rsl] \subseteq \al J.$ and
since $\al J.$ is minimal one has equality,  
$[\rsl \al K. (S_0) \rsl] = \al J.$. The same reasoning applies to 
all larger finite dimensional and non--degenerate subspaces $S \supset S_0$.
But, as we shall see, one has the strict inclusion
$[\rsl \al K. (S) \rsl] \subsetneq [\rsl \al K. (S_0) \rsl]$
whenever $S_0 \subsetneq S$. Hence no nonzero ideal of $\rsl$ can  
be minimal if $X$ is infinite dimensional.   

For the proof of the above inclusion we recall that   
$\al K.(S)$ is contained in all nonzero ideals of $\al R.(S,\sigma)$,  
so one clearly has 
$\al K.(S) \subseteq [\al R.(S,\sigma) \, \al K.(S_0) \, \al
R.(S,\sigma)] \subseteq [\rsl \, \al K.(S_0) \, \rsl]$ and
consequently $[\rsl \al K. (S) \rsl] \subseteq [\rsl \al K. (S_0)
\rsl]$ if $S_0 \subseteq S$.  In order to see that this 
inclusion is strict if $S_0 \subsetneq S$ we pick any pure
state $\omega_S$ on $\al R.(S,\sigma)$ which is regular on 
$\al R.(S_0,\sigma)$ and has in its kernel the resolvents 
$R(\lambda,f)$ with $\lambda \in \R \backslash \{0\}$, 
$f \in S \backslash S_0$,
cf.~Proposition~\ref{PrimIdeals}. We extend $\omega_S$  to a  
state $\omega$ on $\rsl$ and consider its GNS--representation 
$(\pi_\omega, \al H._\omega)$. According to Lemma \ref{SingularStates}(ii) 
the resolvents $R(\lambda,f)$ and therefore also the spaces 
$[R(S,\sigma) \, R(\lambda,f) \, R(S,\sigma)] $,
$f \in S \backslash S_0$ are contained in the kernel of $\pi_\omega$.
But each space $[R(S,\sigma) \, R(\lambda,f) \, R(S,\sigma)] $ is a nonzero
ideal of $R(S,\sigma)$ and consequently $\al K.(S)$ lies in the
kernel of $\pi_\omega$ as well. So 
$\pi_\omega \rest [\rsl \, \al K. (S) \, \rsl] = 0$. On the other hand 
$\omega \rest \al  K.(S_0) = \omega_S \rest \al  K.(S_0)$ is regular by 
construction, hence $\pi_\omega \rest [\rsl \al K. (S_0) \rsl]$ 
is different from $0$. So the respective ideals are different,
which completes the proof of the statement.
\end{beweis}

\vspace*{-4mm}
Having clarified the structure of the basic principal ideals, of their 
intersections and of the minimal ideals 
of the resolvent algebras $\rsl$ we will determine now the   
maximal ideals. We will show that these are  
generated by certain specific subfamilies of the 
basic principal ideals. 
\begin{teo} (i) Let $Z \subset X$ be an isotropic subspace, \ie
$\sigma(f,g) = 0$ for all $f,g \in Z$, let $\al A.(Z) \subset \rsl$ be 
the abelian C*--algebra generated by $R(\lambda,f)$ where 
$\lambda \in \R \backslash \{0\}$, \mbox{$f \in Z$}
and let $\chi$ be a pure state on $\al A.(Z)$ which has none of 
these generating resolvents in its kernel. The 
closed two--sided ideal $\al J.$ generated by 
$\{(R(\lambda,f) - \chi
(R(\lambda,f)) \, \un) : \lambda \in \R \backslash \{0\}, \,
f \in Z \}$ and $\{ R(\mu,g) : \mu \in \R \backslash \{0\}, \,
g \in X \backslash Z \}$ is a proper maximal ideal of $\rsl$. 

(ii) Conversely, let $\al J. \subset \rsl$ be a proper 
maximal ideal of $\rsl$. 
There exists an  isotropic subspace $Z \subset X$ and a pure state 
$\chi$ on the abelian C*--algebra $\al A.(Z)$ 
generated by $R(\lambda,f)$ with   
$\lambda \in \R \backslash \{0\}$, $f \in Z$, which has none of these
generating resolvents in its kernel, such that 
$\al J.$ coincides with the corresponding ideal defined in (i).
\end{teo}

\pagebreak
\begin{beweis}
(i) For the proof that $\al J.$ is a proper ideal we construct 
a representation  $\pi_{\al J.}$ of $\rsl$ which has $\al J.$ in its kernel.
This representation 
acts on the one--dimensional Hilbert space $\C$ and is fixed 
by setting $\pi_{\al J.}(R(\lambda,f)) \doteq \chi(R(\lambda,f)) $
for $\lambda \in \R \backslash \{0\}$, $f \in Z$ and 
$\pi_{\al J.}(R(\mu,g)) \doteq 0$ for $ \mu \in \R \backslash \{0\}, \,
g \in X \backslash Z $. Since $\chi : \al A.(Z) \rightarrow \C$
is a homomorphism one easily verifies that  
$\pi_{\al J.}$ extends by linearity and multiplicativity
to the *--algebra $\al R._0$ generated by the resolvents, \ie
its definition is consistent with the defining relations (\ref{Resolv})
to (\ref{Riden}). By continuity it can therefore be extended to
all of $\rsl$. By construction
$\pi_{\al J.} \rest \al J. = \{0\}$, \ie $\al J.$ is a proper
ideal. For the proof that it
is maximal we recall that $\rsl$ is the closure of 
all polynomials formed out of resolvents and the identity operator. 
Replacing in any given polynomial the resolvents 
$R(\lambda,f)$ by 
$\big((R(\lambda,f) - \chi(R(\lambda,f)) \, \un) + 
\chi(R(\lambda,f)) \, \un \big)$
if $\lambda \in \R \backslash \{0\}$, 
$f \in Z$   and keeping the other resolvents 
untouched it is apparent that
$[\al J. + \C \, \un] = \rsl$. Hence the codimension of  
the proper ideal $\al J.$ is $1$, so it 
cannot be extended any further to a proper ideal, \ie it 
is maximal. 
\\
(ii) For the proof of the second statement 
we make use of the fact that every maximal ideal $\al J.$ is
primitive, \ie kernel of some irreducible representation 
$(\pi_{\al J.},\al H._{\al J.})$, cf.~\cite[Subsec.~II.6.5.3]{Bla}. 
It is an immediate consequence 
of Lemma \ref{Xdecomp}(iii) and the spectral theorem for 
normal operators that for each resolvent $R(\lambda,f) \in \rsl$ with 
$\lambda \in \R \backslash \{0\}$, 
$f \in X$ there is a spectral value $\rho \in \spec(R(\lambda,f))$
such that the inverse of $\pi_{\al J.}(R(\lambda,f) - \rho \,
\un)$ does not exist as a bounded operator on $\al H._{\al J.}$. 
As a matter of fact, for any such $\rho$ one must have 
$(R(\lambda,f) - \rho \, \un) \in \al J.$. For otherwise there exist  
in view of Proposition \ref{PrincId}(ii) and the 
fact that $\al J.$ is a maximal ideal 
sequences $\{J_n \in \al J.\}_{n \in \N}$, $\{R_n \in \rsl\}_{n \in \N}$
such that $J_n + (R(\lambda,f) - \rho \, \un) \, R_n \rightarrow \un$
in the norm topology. Since the invertible elements in a 
C*--algebra form an open set it follows that for sufficiently large 
$n$ one has $(J_n + (R(\lambda,f) - \rho \, \un) \, R_n)^{-1} \in
\rsl$. Thus there are $J \in \al J.$, $R \in \rsl$
such that $J + (R(\lambda,f) - \rho \, \un) \, R = \un$. 
But this implies 
$\pi_{\al J.}(R(\lambda,f) - \rho \, \un) \, \pi_{\al J.}(R) =
\un_{\al H._{\al J.}}$, \ie $ \pi_{\al J.}(R(\lambda,f) - \rho \,
\un)$ has a bounded right inverse. In a similar manner one shows that
$ \pi_{\al J.}(R(\lambda,f) - \rho \, \un)$ has also a bounded left
inverse, which is in conflict with the initial choice of $\rho$.
Hence $(R(\lambda,f) - \rho \, \un) \in \al J.$,
the kernel of $\pi_{\al J.}$, and consequently
$\pi_{\al J.}(R(\lambda,f)) \in \C$ for 
$\lambda \in \R \backslash \{0\}$, $f \in X$, \ie
$(\pi_{\al J.}, \al H._{\al J.} )$ is a one--dimensional 
representation of $\rsl$. 

According to Lemma \ref{Xdecomp}(ii) the set $Z \subset X$ for which 
$\pi_{\al J.}(R(\lambda,f)) \in \C \backslash \{0\}$ if 
$\lambda \in \R \backslash \{0\}$, $f \in Z$, is an isotropic
subspace. Let $\al A. (Z) \subset \rsl$ be the abelian C*--algebra
generated by the corresponding resolvents and let 
$\chi \doteq \pi_{\al J.} \rest \al A. (Z)$.
Since \mbox{$\pi_{\al J.} : \al A. (Z) \rightarrow \C$} is 
a homomorphism, $\chi$ is a pure state which 
does not vanish on any of the generating \mbox{resolvents}. 
According to the preceding step the ideal generated by 
 $\{ (R(\lambda,f) - \chi
(R(\lambda,f)) \, \un) : \lambda \in \R \backslash \{0\}, \,
f \in Z \}$ and $\{ R(\mu,g)  : \mu \in \R \backslash \{0\}, \,
g \in X \backslash Z \}$ is a proper maximal ideal of $\rsl$. 
But, as has been shown, these generating elements are also 
contained in the given 
maximal ideal $\al J.$. So the two ideals must coincide, completing
the proof of the statement. \end{beweis}

As has been shown in the preceding proof, each maximal ideal 
of $\rsl$ coincides with the kernel of a one--dimensional representation and
hence corresponds to a character. Moreover, 
these characters separate the generating resolvents
$\{ \rlf : \lambda \in \Ro, \, f \in X \}$, so there are many of them.
We conclude this analysis of the ideal structure of the 
resolvent algebra with a brief discussion of its commutator ideal,
\ie the ideal which is generated by the commutators of all of its elements. 

\vspace*{-2mm} 
\begin{pro}
The ideal $\al J._c$ generated by 
$\{ [R, R^\prime] : R, R^\prime \in \rsl \}$ is proper. 
It coincides with the ideal generated by 
$\{ R(\lambda,f) R(\mu,g) :  
\lambda, \mu \in \R \backslash \{0\} , \ f,g \in X \ 
\mbox{with} \ \sigma(f,g) \neq 0 \}$ and is contained in all
maximal ideals of $\rsl$. 
\end{pro}
\begin{beweis}
Since the algebra $\al R._0$ of all polynomials in the basic resolvents 
is dense in $\rsl$ it follows that $\al J._c$ coincides with
the ideal generated by 
\mbox{$\{ [R(\lambda,f), R(\mu,g)] :
\lambda, \mu \in \R \backslash \{0\} , \ f,g \in X \}$}.  
According to relation (\ref{Rccr}) 
the latter ideal in turn coincides with the ideal generated by 
\mbox{$\{ R(\lambda,f) R(\mu,g)^2  R(\lambda,f) : 
\lambda, \mu \in \R \backslash \{0\}, 
\ f,g \in X \ \mbox{with} \ \sigma(f,g) \neq 0 \}$}. Now 
by Proposition~\ref{IdIntersect} the ideals generated by 
$R(\lambda,f) R(\mu,g)^2  R(\lambda,f)$, being equal to 
$\al I.(\lambda,f;0) \cap \al I.(\mu,g;0) $, 
coincide with the ideals generated by $R(\lambda,f) R(\mu,g)$
for $\lambda, \mu \in \Ro$ and $f,g \in X$.
Hence $\al J._C$ is equal to the ideal generated by 
$\{ R(\lambda,f) R(\mu,g) : 
\lambda, \mu \in \R \backslash \{0\}, 
\ f,g \in X \ \mbox{with} \ \sigma(f,g) \neq 0 \}$.

It remains to show that the ideal $\al J._c$ 
is contained in all maximal ideals determined in 
the preceding theorem (which also implies that $\al J._c$ is proper). 
But this instantly follows from the fact, established above, that 
all maximal ideals coincide with the kernels of 
one--dimensional representations of $\rsl$. These 
annihilate all commutators and consequently also $\al J._c$. 
\end{beweis}

Denoting by $R_c(\lambda,f)$ the class of 
$R(\lambda,f)$ modulo $\al J._c$ where $\lambda \in \R \backslash
\{0\}$, $f \in X$ one easily checks that these 
operators satisfy the defining relations 
(\ref{Resolv}) to (\ref{Riden}) with $\sigma \equiv 0$ and they 
generate the abelian C*--algebra $\rsl /\al J._c$. The
symplectic form $\sigma$ only enters in the additional relation
$R_c(\lambda,f) R_c(\mu,g) = 0$ if $\sigma(f,g) \neq 0$.  
Hence the operators assigned to such incompatible elements of $X$ 
have disjoint spectral supports, reflecting the incommensurability 
of the underlying quantum observables in the abelian quotient
algebra.

\section{Conclusions}
\setcounter{equation}{0}
\label{conclusions}

In the present investigation we have clarified the ideal structure
of the resolvent algebra $\rsl$. All of its ideals 
are built in a simple, physically significant  
manner from principal ideals generated by 
the basic resolvents. Moreover, the nesting of its primitive ideals 
encodes precise information about the dimension of the symplectic
space $X$, \ie of the size of the underlying quantum system. 
There is a sharp algebraic 
distinction between quantum systems with finitely many particles, where 
$\dim(X) < \infty$ and the resolvent algebra is postliminal (type $I$),
and the case of quantum field theory, respectively 
infinitely many particles, where $\dim(X) = \infty$ and the resolvent algebra
is no longer postliminal, but it is still nuclear.

Another prominent difference between these two cases consists of the 
following fact: If  $\dim(X) < \infty$ the resolvent algebra 
$\rsl$ contains a non--trivial minimal ideal $\al K.$ which is isomorphic
to the compacts. This ideal carries the regular representations 
of $\rsl$ in the sense that these are precisely the unique extensions
of the (nondegenerate) representations of~$\al K.$. If 
$\dim(X) = \infty$ there exists no such  ideal in $\rsl$, however. 
Yet since $\rsl$ is the C*--inductive limit of its subalgebras 
$\al R.(S, \sigma)$ for all finite dimensional nondegenerate 
subspaces $S \subset X$ there is a certain substitute. 
Each subalgebra $\al R.(S, \sigma) \subset \rsl$ contains its own 
compact ideal $\al K.(S)$. 
Let $\al L. \subset \rsl$ be the C*--algebra  
which is generated by all algebras $\al K.(S)$, $S \subset X$. 
One can show that $\al L.$ is a bimodule for $\rsl$ and 
lies in the kernel of the one--dimensional
representation which annihilates all resolvents, hence
$\al L.$ is a proper ideal. It carries the regular representations 
of $\rsl$ in the sense that these are precisely the unique extensions
of the representations of~$\al L.$ whose restrictions to 
all subalgebras $\al K.(S) \subset \al L.$, $S \subset X$ 
are nondegenerate. 

The properties of the ideal~$\al L.$ and of
its subalgebras are also a key to the 
construction of interesting automorphism groups (dynamics)
of the resolvent algebra \cite{BuGr2}. We hope to return to this 
physically important issue in a future publication.

\section*{Acknowledgments}

I am grateful to Hendrik Grundling for numerous
discussions and valuable contributions which helped to
clarify this topic. I would also like to thank Ralf Meyer
for a crucial  remark on invariants and ideals of
C*--algebras.


\bigskip

\begin{thebibliography}{10}
\small 

\bibitem{Ace}
F.\ Acerbi, G.\ Morchio and F.\ Strocchi:
Nonregular representations of CCR algebras and algebraic fermion bosonization.
Proceedings of the XXV Symposium on Mathematical Physics (Toru\'n, 1992).
Rep. Math. Phys.~{\bf 33} no. 1-2, 7--19 (1993).

\bibitem{Bla}
R.\ Blackadar: \textit{Operator Algebras}. Springer 2006

\bibitem{BuGr}
D.\ Buchholz and H.\ Grundling: Algebraic supersymmetry: A case study.
Commun. Math. Phys.~{\bf 272},  699--750 (2007)

\bibitem{BuGr2}
D.\ Buchholz and H.\ Grundling: The resolvent algebra: A new 
approach to canonical quantum systems.
Journal of Functional Analysis~{\bf 254},  2725--2779 (2008)

\bibitem{BuGr3}
D.\ Buchholz and H.\ Grundling:  
Lie algebras of derivations and resolvent algebras. 
Comm. Math. Phys.~{\bf 320}, 455--467 (2013)  

\bibitem{Cost}
R.\ Costello: The mathematics of the BRST-constraint method.
\ e--print \ arXiv:0905.3570

\bibitem{Dix}
J.\ Dixmier: \textit{C*--algebras}. 
North Holland Publishing Co.  1977

\bibitem{FaVB}
M.\ Fannes and A.\ Verbeure: On the time evolution automorphisms of the
CCR--algebra for quantum mechanics.
Commun. Math. Phys.~{\bf 35}, 257--264 (1974)

\bibitem{Gr3}
H.\ Grundling: A group algebra for inductive limit groups. Continuity problems
  of the canonical commutation relations. Acta Applicandae 
Mathematicae~\textbf{46}, 107--145 (1997)

\bibitem{GrHu3}
H.\ Grundling and C.A.\ Hurst:
 A note on regular states and supplementary conditions. Lett.
  Math. Phys.~\textbf{15}, 205--212 (1988) [Errata: ibid. {\bf 17},
  173--174 (1989)]

\bibitem{Man}
J.\ Manuceau: C*-alg\`ebre de relations de commutation. Annales de l'Institut
H. Poincar\'e (A)~{\bf 8}, 139--161 (1968)

\bibitem{Ped}
Pedersen, G.K.: \textit{C${}^{\, *}$--Algebras and their Automorphism Groups}.
  Academic Press 1989

\bibitem{Ro}
P.L.\ Robinson: \ Symplectic pathology. \
Q.\ J.\ Math.~{\bf 44}, 101--107 (1993) 

 \bibitem{Yos}
 K.\ Yosida: \textit{Functional Analysis}. 
Springer  1980
\end{thebibliography}
\end{document}